\documentclass[11pt,reqno]{elsarticle}

\makeatletter
\def\ps@pprintTitle{%
 \let\@oddhead\@empty
 \let\@evenhead\@empty
 \def\@oddfoot{}%
 \let\@evenfoot\@oddfoot}
\makeatother

\usepackage{amsmath,amsthm,amscd,amsfonts,amssymb,mathtools}
\usepackage{mathabx,breqn}
\usepackage{booktabs,cool}
\usepackage{parskip}
\usepackage[top=1.1in, bottom=1.2in, left=1in, right=1in]{geometry}

\DeclareMathOperator{\cotan}{cotan}
\DeclareMathOperator{\real}{Re}
\DeclareMathOperator{\imag}{Im}

\newtheorem{theorem}{Theorem}

\providecommand{\e}[1]{\ensuremath{\times 10^{#1}}}

\makeatletter
\g@addto@macro\normalsize{%
  \setlength\abovedisplayskip{.4em}
  \setlength\belowdisplayskip{.4em}
  \setlength\abovedisplayshortskip{.4em}
  \setlength\belowdisplayshortskip{.4em}
}

\begin{document}

\begin{frontmatter}
\begin{abstract}
We express a certain complex-valued solution 
of Legendre's differential equation  as the product of 
an oscillatory exponential function and an integral
involving only nonoscillatory elementary functions.
By calculating the logarithmic derivative of this solution, we show
that Legendre's differential equation admits a nonoscillatory phase function.
Moreover, we derive from our expression
an asymptotic expansion useful for evaluating
 Legendre functions of the first and second kinds of large orders,
as well as the derivative of the nonoscillatory phase function.
Numerical experiments demonstrating the properties of our asymptotic expansion are presented.
\end{abstract}


\title
{
On the nonoscillatory phase function for Legendre's differential equation
}

\author[jb]{James Bremer\corref{cor1}}
\ead{bremer@math.ucdavis.edu}
\author[vr]{Vladimir Rokhlin}

\cortext[cor1]{Corresponding author}

\address[vr]{Department of Computer Science, Yale University}
\address[jb]{Department of Mathematics, University of California, Davis}


\end{frontmatter}

\begin{section}{Introduction}

 The Legendre functions  of degree $\nu$
--- that is, the solutions of the
second order  linear ordinary differential equation
\begin{equation}
\psi''(z) - \frac{2 z}{1-z^2} \psi'(z) + \frac{\nu (\nu +1)}{1-z^2} \psi(z) = 0 
\label{introduction:legendre}
\end{equation}
---  
appear in numerous contexts in  physics and applied mathematics.
For instance,
they arise when certain partial differential equations are solved
via separation of variables,  they are often used to represent smooth functions defined on 
bounded intervals, and their roots are the nodes of Gauss-Legendre quadrature rules.

In this article, we give a constructive proof of the existence of a nonoscillatory
phase function for Legendre's differential equation.  
A smooth function $\alpha:(a,b) \to \mathbb{R}$ is a phase
function for the second order linear ordinary differential equation
\begin{equation}
y''(x) + q(x) y(x) = 0 \ \ \mbox{for all} \ \ a < x < b,
\label{introduction:ode}
\end{equation}
where $q$ is smooth and positive on $(a,b)$ (it can have zeros or singularities
at the endpoints of the interval),
if $\alpha'$ does not vanish on $(a,b)$ and the functions 
\begin{equation}
\frac{\cos(\alpha(x))}{\sqrt{\left|\alpha'(x)\right|}}
\label{introduction:u0}
\end{equation}
and
\begin{equation}
\frac{\sin(\alpha(x))}{\sqrt{\left|\alpha'(x)\right|}}
\label{introduction:v0}
\end{equation}
form a basis in the space of solutions of (\ref{introduction:ode}).
When $q$ is positive the solutions of (\ref{introduction:ode})
are oscillatory, whereas when $q$ is negative the solutions of (\ref{introduction:ode})
are roughly increasing and decreasing exponentials.
In the latter case,  (\ref{introduction:u0}) and (\ref{introduction:v0}) 
are not appropriate mechanisms for representing them, hence the requirement that $q$ 
is positive.    The ostensibly more general second order linear ordinary
differential equation
\begin{equation}
\eta''(x) + p_0(x) \eta'(x) + q_0(x) \eta(x)  = 0\ \ \mbox{for all}\ \ a < x < b,
\label{introduction:ode2}
\end{equation}
where the coefficients $p_0$, $q_0$  are real-valued and smooth on $(a,b)$,
can be  reduced to the form (\ref{introduction:ode}) 
with the  function $q$ given by 
\begin{equation}
q(x) =   q_0(x) - \frac{1}{4} (p_0(x))^2 - \frac{1}{2} p_0'(x).
\label{introduction:q}
\end{equation}
This is accomplished  by letting  
\begin{equation}
y(x) = \frac{\eta(x)}{\sqrt{W(x)}},
\end{equation}
%
where  $W$  is defined via
\begin{equation}
W(x) = \exp\left(-\int_{a}^x p_0(u)\ du\right).
\end{equation}
%
%
Assuming that the function $q$ defined by  (\ref{introduction:q}) is positive,
we say $\alpha$ is a phase function for (\ref{introduction:ode2})
if it is a phase function for the transformed
equation (\ref{introduction:ode}).  In this case, the functions
$u$, $v$ defined via the formulas
\begin{equation}
u(x) = \sqrt{\frac{W(x)}{\left|\alpha'(x)\right|}}\ \cos(\alpha(x))
\label{introduction:u2}
\end{equation}
and
\begin{equation}
v(x) =  \sqrt{\frac{W(x)}{\left|\alpha'(x)\right|}}\ \sin(\alpha(x))
\label{introduction:v2}
\end{equation}
form a basis in the space of solutions of the original equation (\ref{introduction:ode2}).
We note that according to Abel's identity, the Wronskian of any 
pair of solutions of (\ref{introduction:ode2}) is a constant multiple
of the function $W$, and that the Wronskian of the pair
defined by (\ref{introduction:u2}) and (\ref{introduction:v2})
is $W$.

It has long been known that certain second order differential equations admit
nonoscillatory phase functions, and that a large class of such
equations do  so in an asymptotic sense. Indeed, under mild conditions on $q$,  
the second order equation
\begin{equation}
y''(x) + \lambda^2 q(x) y(x) = 0\ \ \mbox{for all}\ \ a < x < b
\label{introduction:ode3}
\end{equation}
admits a basis $\{u_1,v_1\}$ of solutions such that
\begin{equation}
u_1(x) 
= \frac{1}{\sqrt{\lambda}\ (q(x))^{1/4}} \cos\left(\lambda \int_{a}^x \sqrt{q(u)}\ du\right)
+ O\left(\frac{1}{\lambda}\right)
\label{introduction:u1}
\end{equation}
and
\begin{equation}
v_1(x) 
 = \frac{1}{\sqrt{\lambda}\ (q(x))^{1/4}} \sin\left(\lambda \int_a^x \sqrt{q(u)}\ du\right)
+ O\left(\frac{1}{\lambda}\right).
\label{introduction:v1}
\end{equation}
This standard result can be found in \cite{Olver} and \cite{Fedoryuk}, among many other sources.
Assuming that the coefficient $q$ is nonoscillatory,  the phase function
\begin{equation}
\alpha_1(x) = \lambda \int_a^x \sqrt{q(u)}\ du
\end{equation}
%
associated with (\ref{introduction:u1}) and (\ref{introduction:v1})
obviously is as well. 
Expressions of this type are generally known as WKB approximations,
and higher order analogues of them can be constructed.  
The resulting formulas, however, are quite unwieldy --- 
indeed, the $n^{th}$ order WKB expansion,
which achieves $O\left(\lambda^{-n}\right)$ accuracy,
involves complicated combinations of 
the derivatives of $q$ of orders $0$ through $n-1$ and various
noninteger powers
of $q$.
Perhaps more importantly,  efficient and accurate numerical methods for their
computation are lacking.

It is shown in \cite{Heitman-Bremer-Rokhlin,Bremer-Rokhlin} that
under mild conditions on $q$, 
 there exist a positive constant $\mu$,
 a function $\alpha_\infty$ which is roughly as oscillatory
as $q$, and a basis $\{u_\infty,v_\infty\}$
in the space of solutions of the differential equation
(\ref{introduction:ode3}) such that 
\begin{equation}
u_\infty(x) =
\frac{\cos(\alpha_\infty(x))}{\sqrt{\left|\alpha_\infty'(x)\right|}}  + O\left(\exp(-\mu \lambda)\right)
\end{equation}
and
\begin{equation}
v_\infty(x) =
\frac{\sin(\alpha_\infty(x))}{\sqrt{\left|\alpha_\infty'(x)\right|}}  + O\left(\exp(-\mu \lambda)\right).
\end{equation}
The function $\alpha_\infty$ can be represented using
various series expansions (such as piecewise polynomial expansions) the number
of terms of  which do not depend on the parameter $\lambda$.  In order words,
nonoscillatory phase functions represent the solutions of equations 
of the form (\ref{introduction:ode3}) with  $O(\exp(-\mu\lambda)$ accuracy 
using $O(1)$-term expansions.
Moreover,  a reliable and efficient
numerical algorithm for the computation of $\alpha_\infty$ which runs in time
independent of $\lambda$ and only requires knowledge of the values
of $q$ on the interval $(a,b)$ is introduced in \cite{BremerKummer}.
Much like standard WKB estimates,
the results of \cite{Heitman-Bremer-Rokhlin,Bremer-Rokhlin,BremerKummer}
easily apply to large class of differential equations whose coefficients
are allowed to vary with $\lambda$ --- that is, equations of the more
general form
\begin{equation}
y''(x) + q(x,\lambda) y(x) = 0
\end{equation}
--- assuming that $q$ satisfy certain innocuous conditions independent of $\lambda$.

The framework of \cite{Heitman-Bremer-Rokhlin,Bremer-Rokhlin,BremerKummer}
applies to (\ref{introduction:legendre})
and it can be used to, among other things, evaluate Legendre functions of large
orders and their zeros in time independent of degree.  However, in cases
like Legendre's differential equation where an exact nonoscillatory phase function
exists, it is advantageous to derive formulas for it which are as explicit 
as possible.    
To that end, in this article, we derive an integral representation of
a particular solution $\psi_\nu$ of Legendre's differential equation.  
This is done in Section~\ref{section:integral}, after certain
preliminaries are dispensed with in Section~\ref{section:preliminaries}.
In Section~\ref{section:nonoscillatory}, our integral representation
of $\psi_\nu$ is used to prove
the existence of a nonoscillatory phase function $\alpha_\nu$
for Legendre's differential equation.
In Section~\ref{section:asymptotics}, we 
 derive an asymptotic formula for Legendre functions of large degrees,
and for the derivative of the nonoscillatory phase function $\alpha_\nu$.
Numerical experiments demonstrating the properties of these expansions are
discussed in Section~\ref{section:numerics}.  We conclude this article
with brief remarks in Section~\ref{section:conclusion}.

\end{section}


\begin{section}{Preliminaries}

\begin{subsection}{An elementary observation regarding phase functions}
\label{section:preliminaries:phase}

Suppose that $u$, $v$ form a basis in the space of solutions
of (\ref{introduction:ode2}).  Then, since $u$ and $v$ 
cannot simultaneously vanish,
 (\ref{introduction:u2}) and (\ref{introduction:v2})
hold for all $x \in (a,b)$ if and only if 
\begin{equation}
\tan(\alpha(x)) = \frac{v(x)}{u(x)}
\label{preliminaries:1}
\end{equation}
for all $x \in (a,b)$ such that $u(x) \neq 0$ and 
\begin{equation}
\cotan(\alpha(x)) = \frac{u(x)}{v(x)}
\label{preliminaries:2}
\end{equation}
for all $x \in (a,b)$ such that $v(x) \neq 0$.
By differentiating 
(\ref{preliminaries:1}) and (\ref{preliminaries:2}) we find that
(\ref{introduction:u2}) and (\ref{introduction:v2}) are satisfied
if and only if 
\begin{equation}
\alpha'(x) = \frac{u(x)v'(x) - u'(x)v(x)}{(u(x))^2 + (v(x))^2}
\label{preliminaries:alphap}
\end{equation}
for all $a < x < b$.  In particular, if $u$, $v$ form a basis in the space
of solutions of (\ref{introduction:ode2}), then any function $\alpha$ whose
derivative satisfies
(\ref{preliminaries:alphap}) is a phase function for (\ref{introduction:ode2}).

We will frequently make use of the following elementary theorem:
\vskip 1em
\begin{theorem}
If $u$, $v$ are real-valued solutions of (\ref{introduction:ode2}) and 
\begin{equation}
\psi(x) = u (x) + i v(x),
\end{equation}
then the imaginary part of the logarithmic derivative $\psi'(x) / \psi(x)$
of $\psi$ is the derivative of a phase function for (\ref{introduction:ode2})
\label{preliminaries:theorem1}
\end{theorem}
It is easily established by observing that
\begin{equation}
\frac{\psi'(x)}{\psi(x)} = 
\frac{u(x) u'(x) + v(x) v'(x)}{(u(x))^2 + (v(x))^2}
+ i
\frac{u(x) v'(x) - u'(x) v(x)}{(u(x))^2 + (v(x))^2}
\label{preliminaries:alphap0}
\end{equation}
and comparing (\ref{preliminaries:alphap}) and (\ref{preliminaries:alphap0}).

\end{subsection}


\begin{subsection}{Gauss' hypergeometric function}

For any complex-valued parameters $a$, $b$, and $c$ such that $c \neq 0,-1,-2,\ldots$ 
we use $\Hypergeometric{2}{1}{a,b}{c}{z}$ to denote Gauss' hypergeometric function,
which is defined by the formula
\begin{equation}
\Hypergeometric{2}{1}{a,b}{c}{z}
 = \sum_{n=0}^\infty \frac{\left(a\right)_n \left(b\right)_n}{\left(c\right)_n} \frac{z^n}{\Gamma(n+1)}.
\label{preliminaries:hypergeometric}
\end{equation}
Here, $(z)_n$ is the Pochhammer symbol
\begin{equation}
z (z+1) (z+2) \cdots (z+n-1).
\end{equation}
Since the series (\ref{preliminaries:hypergeometric})
converges absolutely for all $|z|<1$, $\Hypergeometric{2}{1}{a,b}{c}{z}$ 
is an analytic function  on the unit disk in the complex plane.

\end{subsection}

\begin{subsection}{Legendre functions}

We refer to the function
\begin{equation}
P_\nu(z) = 
\Hypergeometric{2}{1}{-\nu,\nu+1}{1}{\frac{1-z}{2}}
\label{definitions:pnu}
\end{equation}
obtained by employing the Frobenius (series solution) method to construct
a solution of (\ref{introduction:legendre}) which is analytic in a neighborhood
of the regular singular point $z=1$  as the Legendre function of the first kind of degree $\nu$.  
When $\nu$ is an integer,
the series (\ref{definitions:pnu}) terminates and $P_\nu$ is entire; otherwise,
it can be analytically continued to the cut plane $\mathbb{C}\setminus\left(-\infty,-1\right]$
(see, for instance,
Section~2.10 of \cite{HTFI}).
Similarly, we refer to the solution
\begin{equation}
Q_\nu(z) = 
(2z)^{-\nu-1} \sqrt{\pi} \frac{\Gamma(\nu+1)}{\Gamma(\nu+3/2)}
\ \Hypergeometric{2}{1}{\frac{\nu}{2}+1,\frac{\nu}{2}+\frac{1}{2}}{\nu+\frac{3}{2}}{\frac{1}{z^2}}
\label{deinfitions:qnu}
\end{equation}
%
of (\ref{introduction:legendre}) 
 as the Legendre function of the second kind of degree $\nu$.   It can be
analytically continued 
to the cut plane $\mathbb{C}\setminus\left(-\infty,1\right]$.  
Following the standard convention  (see, for example, \cite{NIST:DLMF,HTFI,Gradshteyn}), we define
$Q_\nu(x)$ for $x \in (-1,1)$  via the formula
\begin{equation}
Q_\nu(x) = \lim_{y \to 0^+} \frac{Q_\nu(x+iy) + Q_\nu(x-iy)}{2}.
\end{equation}

It can be verified readily that while $Q_\nu$ is not an element of any
of the classical Hardy spaces on the upper half of the complex plane $\mathbb{H}$, 
it is an element of the space $H^+$ of functions which are analytic on $\mathbb{H}$
and uniformly bounded by a polynomial on every half plane of the form $\imag(z) \geq y_0 > 0$.
That is, a function $F$ analytic on $\mathbb{H}$ is an element of $H^+$
if for every $y_0 > 0$ there exists a polynomial $p_{y_0}$ such that
\begin{equation}
\sup_{y \geq y_0} \left| F(x+iy) \right| < p_{y_0}(x)
\ \ \mbox{for all} \ \ x\in\mathbb{R}.
\end{equation}
The space $H^+$ generalizes the classical Hardy spaces $H^p$ and many of the useful
properties of Hardy spaces still apply in  this setting.
For instance, if $1 \leq p <\infty$ then the
elements of $H^p$  are the Fourier transforms
of $L^p$ functions supported on the half line $[0,\infty)$, while
those of $H^+$ are the Fourier transforms of distributions supported
on $[0,\infty)$.
 The classical theory of Hardy spaces
is discussed in the well-known textbook \cite{Koosis} 
and the properties of $H^+$ are described in some detail in \cite{Beltrami-Wohlers}.

   According to Formula~(9) in Section~3.4 of \cite{HTFI} 
(see also 8.732.5 in \cite{Gradshteyn}),
\begin{equation}
\lim_{y\to 0^+} Q_\nu(x+iy) = \varphi_\nu(x),
\label{definitions:fv}
\end{equation}
where
\begin{equation}
\varphi_\nu(x) = Q_\nu(x) - i \frac{\pi}{2} P_\nu(x)
\label{preliminaries:phi}
\end{equation}
for all $-1 < x <1$.  
Using the fact that $Q_\nu \in H^+$ and the observation
that it  has no zeros in the upper half of the complex plane
(this is established, for instance, in \cite{Runkel}  via the argument principle),
it can be shown that the logarithmic derivative of the function
$\varphi_\nu$  is  nonoscillatory. 
It then follows from Theorem~\ref{preliminaries:theorem1}
 that Legendre's differential
equation admits a nonoscillatory phase function.  The details of this argument
are beyond the scope of this article, but we mention it because
of its strong relation to the approach  taken here.

\label{section:preliminaries:legendre}
\end{subsection}

\begin{subsection}{A transformation of Legendre's differential equation}

It is easy to verify that
the functions  $P_\nu(\cos(\theta))$
and $Q_\nu(\cos(\theta))$ satisfy the second order differential equation
\begin{equation}
\psi''(\theta) + \cot(\theta) \psi'(\theta) + \nu(\nu+1) \psi(\theta) = 0
\label{preliminaries:legendre}
\end{equation}
on the interval $(0,\pi/2)$.
We will, by a slight abuse of terminology, refer to both (\ref{introduction:legendre})
and (\ref{preliminaries:legendre}) as Legendre's differential equation.
In all cases, though, it will be clear from the context which specific form of Legendre's
differential equation is intended.

\end{subsection}

\begin{subsection}{Bessel functions}


For nonnegative integers $n$, 
 the Bessel function of the first kind of order $n$ is the entire function
defined via the series
\begin{equation}
J_n(z) = \sum_{j=0}^\infty \ \frac{(-1)^j}{\Gamma(j+1)^2} \left(\frac{z}{2}\right)^{2j+n},
\label{definitions:j0def}
\end{equation}
while the Bessel function of the second kind of order $n$ is 
\begin{dmath}
Y_n(z) = 2 J_n(z) \log\left(\frac{z}{2}\right)
-\sum_{j=0}^{n-1} \frac{\Gamma(n-j)}{\Gamma(j+1)} \left(\frac{z}{2}\right)^{2j-n}\\
- \sum_{j=0}^\infty \frac{1}{\Gamma(j+1) \Gamma(j+n+1)} \left(\frac{z}{2}\right)^{n+2j}
\left(\psi(j+1)+\psi(j+n+1)\right),
\label{preliminaries:y0}
\end{dmath}
where $\psi$ denotes the logarithmic derivative of the gamma function.
%
We use the principal branch of the logarithm in (\ref{preliminaries:y0}),
so that $Y_n$ is defined on the cut plane $\mathbb{C}\setminus\left(-\infty,0\right]$.

The Hankel function of the first kind of order $n$ is
\begin{equation}
H_n(z) = J_n(z) + i Y_n(z)
\label{preliminaries:hankel}
\end{equation}
and we let $K_n$ denote the modified Bessel function of the first kind
of order $n$.  For complex-valued $z$ such that $-\pi < \arg(z) \leq \pi/2$ ,
$K_n$  is defined in terms of  $H_n$ via the formula
\begin{equation}
K_n(z) =  \frac{\pi}{2} \exp\left(n \frac{\pi}{2} i \right)i H_n(iz).
\label{preliminaries:k0}
\end{equation}
It is easily verified that
\begin{equation}
\frac{d}{dz} H_0(z) = -H_1(z).
\label{preliminaries:h0der}
\end{equation}

\end{subsection}


\begin{subsection}{The Lipschitz-Hankel integrals}

The formulas
\begin{equation}
\frac{1}{\Gamma(\nu+1)}\int_0^\infty \exp(- \cos(\theta) x) J_0( \sin(\theta) x) x^\nu\ dx
=  P_\nu(\cos(\theta))
\label{preliminaries:lipschitz1}
\end{equation}
and
\begin{equation}
\frac{1}{\Gamma(\nu+1)}
\int_0^\infty \exp(- \cos(\theta) x) Y_0(\sin(\theta) x) x^\nu\ dx
= - \frac{2}{\pi}  Q_\nu(\cos(\theta)),
\label{preliminaries:lipschitz2}
\end{equation}
hold when $0 < \theta < \pi/2$ and $\real(\nu) >-1$.
Derivations of (\ref{preliminaries:lipschitz1}) and (\ref{preliminaries:lipschitz2})
can be found in Chapter~13 of \cite{Watson}; they can also be
found as  Formulas~6.628(1) and 6.628(2) in \cite{Gradshteyn}.


\end{subsection}

\begin{subsection}{The Laplace transforms of $\left(x^2 + bx\right)^{-1/2}$ and $x^\nu$}

For complex numbers $z$ and $b$ such that $\real(z) > 0$ and 
$|\arg(b)| < \pi$,
\begin{equation}
\int_0^\infty \frac{\exp(-z x)}{\sqrt{x^2+bx}}\ dx = 
\exp\left(\frac{b z}{2}\right)
K_0\left(\frac{bz}{2}\right)
\label{preliminaries:laplace1_0}
\end{equation}
(as usual, we take the principal branch of the square root function). 
A careful derivation of
(\ref{preliminaries:laplace1_0})
 can be found in Section~7.3.4 of \cite{HTFII} and 
it  appears as  Formula~3.383(8) in \cite{Gradshteyn}.
By combining (\ref{preliminaries:laplace1_0})  with (\ref{preliminaries:k0}) we see that 
\begin{equation}
\int_0^\infty \frac{\exp(-z x)}{\sqrt{x^2+bx}}\ dx = 
\frac{\pi}{2} i \exp\left(\frac{b z}{2}\right)
H_0\left(i\frac{bz}{2}\right)
\label{preliminaries:laplace1}
\end{equation}
whenever   $\real(z) > 0$, $-\pi < \arg(bz) \leq \pi/2$ and 
$|\arg(b)|<\pi$.


It is abundantly  well known that for any complex numbers $z$ and $\nu$ such that
$\real(\nu) > -1$ and $\real(z) > 0$,
\begin{equation}
\int_0^\infty \exp(-zx) x^\nu\ dx = \frac{\Gamma(1+\nu)}{z^{1+\nu}};
\label{preliminaries:laplace2}
\end{equation}
this identity  appears, for instance,  as Formula~3.381(4) in \cite{Gradshteyn}.
\end{subsection}

\begin{subsection}{Steiltjes' asymptotic formula for Legendre functions}
For all positive real $\nu$ and all $0 < \theta < \pi/2$, 
\begin{equation}
P_\nu(\cos(\theta)) =
\left(\frac{2}{\pi \sin(\theta)}\right)^{1/2} \ 
\sum_{k=0}^{M-1}
C_{\nu,k}
\frac{\cos\left( \alpha_{\nu,k} \right)}{\sin(\theta)^k}
+ R_{\nu,M}(\theta),
\label{preliminaries:stieltjes}
\end{equation}
where
\begin{equation}
\alpha_{\nu,k} = \left(\nu + k + \frac{1}{2} \right) \theta - 
\left(k+\frac{1}{2}\right)\frac{\pi}{2},
\end{equation}
\begin{equation}
C_{\nu,k} =  \frac{\left(\Gamma\left(k+\frac{1}{2}\right)\right)^2 \Gamma(\nu+1)}
{\pi 2^k \Gamma\left(\nu+k+\frac{3}{2}\right)\Gamma\left(k+1\right)},
\end{equation}
and
\begin{equation}
\left|R_{\nu,M}(\theta)\right| \leq \left(\frac{2}{\pi\sin(\theta)}\right)^{1/2}
\frac{C_{\nu,M}}{\sin(\theta)^M}.
\label{preliminaries:errorbound}
\end{equation}
This approximation was introduced by Stieltjes; a derivation of the 
 error bound  (\ref{preliminaries:errorbound}) can be found in Chapter~8 of  \cite{Szego}, among many
other sources.  In \cite{Bogaert-Michiels-Fostier}, an efficient method for computing
the coefficients $C_{\nu,k}$ is suggested.
The coefficient in the first term is given by
\begin{equation}
C_{\nu,0} = \frac{\Gamma\left(\nu+\frac{1}{2}\right)}{\Gamma\left(\nu+\frac{3}{2}\right)}.
\end{equation}
This ratio of gamma functions can be approximated via a series in powers of
$1/\nu$; however, it can be calculated more efficiently by observing that
the related function
\begin{equation}
\tau(x) = \sqrt{x} \frac{\Gamma\left(x + \frac{1}{4} \right)}{\Gamma\left(x + \frac{3}{4}\right)}
\end{equation}
admits an expansion in powers of $1/x^2$.  In particular,
the $7$-term expansion
\begin{equation}
\begin{aligned}
\tau(x) \approx 
1-&\frac{1}{64 x^2}+\frac{21}{8192 x^4}-\frac{671}{524288 x^6}+\frac{180323}{134217728
   x^8}
\\ 
&\ \ \ \ \
-\frac{20898423}{8589934592 x^{10}}+\frac{7426362705}{1099511627776
   x^{12}}+O\left(\frac{1}{x^{14}}\right)
\end{aligned}
\end{equation}
gives roughly double precision accuracy for all $x > 10$ (see the discussion
in \cite{Bogaert-Michiels-Fostier}).
The coefficient $C_{\nu,0}$ is related to $\tau$ through the formula
\begin{equation}
C_{\nu,0} = \frac{1}{\sqrt{\nu + \frac{3}{4}}} \tau\left(\nu + \frac{3}{4}\right).
\end{equation}
The subsequent coefficients are obtained through the recurrence relation
\begin{equation}
C_{\nu,k+1} = \frac{\left(k+\frac{1}{2}\right)^2}{2 (k+1) \left(\nu+k+\frac{3}{2}\right)} C_{\nu,k}.
\end{equation}

%
%
%
%
\end{subsection}


\label{section:preliminaries}
\end{section}


\begin{section}{An integral representation of a particular solution of Legendre's equation}

For $0 < \theta < \pi/2$ and $\real(\nu) > -1$, we define the function
$\psi_\nu$ via 
\begin{equation}
\psi_\nu(\theta) = P_\nu(\cos(\theta)) - \frac{2}{\pi} i Q_\nu(\cos(\theta)).
\end{equation}
Obviously, $\psi_\nu$ is a solution of Legendre's equation (\ref{preliminaries:legendre})
and it is related to the function $\varphi_\nu$ defined in 
Section~\ref{section:preliminaries:legendre} via the formula
\begin{equation}
\psi_\nu(z) = \frac{2}{\pi} i\  \overline{\varphi_\nu(z)}.
\end{equation}
In this section, we derive
an integral representation of $\psi_\nu$ involving only elementary functions.
We prefer $\psi_\nu$ over $\varphi_\nu$ because it leads to an asymptotic
formula involving the Hankel functions of the first kind of order $0$ instead of 
the Hankel function of the second kind of order $0$.

From (\ref{preliminaries:lipschitz1}), (\ref{preliminaries:lipschitz2}) 
and the definition (\ref{preliminaries:hankel}) of the Hankel function of the first kind of order $0$,
we see  that 
\begin{equation}
\psi_\nu(\theta) =  
\frac{1}{\Gamma(\nu+1)}\int_0^\infty \exp(- \cos(\theta) t) H_0(\sin(\theta) t) t^\nu\ dt
\label{integral:2}
\end{equation}
whenever $ 0 < \theta < \pi/2$ and $\real(\nu) > -1$. We rearrange 
(\ref{integral:2}) as 
\begin{equation}
\psi_\nu(\theta) =  
\frac{1}{\Gamma(\nu+1)}\int_0^\infty 
\exp(- \exp(-i \theta)  t)
\exp(-i\sin(\theta) t)H_0(\sin(\theta) t) 
t^\nu\ dt,
\label{integral:3}
\end{equation}
let $\beta = \exp(i\theta) \sin(\theta)$
and introduce the new variable $w = \exp(-i\theta)t$  to obtain
\begin{equation}
\psi_\nu(\theta) =  
\frac{\exp(i(\nu+1)\theta)}{\Gamma(\nu+1)}\int_0^\infty 
\exp(- w) \exp(-i\beta w)H_0(\beta w)  w^\nu\ \ dw.
\label{integral:4}
\end{equation}
Now by letting $b = 2i $ and $z= \beta w$ in (\ref{preliminaries:laplace1}), we find
that
\begin{equation}
-\frac{2}{\pi}i\int_0^\infty \frac{\exp(-\beta w  x)}{\sqrt{x^2-2ix}}\ dx = 
\exp\left(-i\beta w\right)
H_0\left(\beta w \right).
\label{integral:4.5}
\end{equation}
Formula~(\ref{preliminaries:laplace1}) holds  since $0 < \arg(z) < \pi/2$  and  $\arg(b) = -\pi/2$.
Inserting (\ref{integral:4.5})  into (\ref{integral:4})  yields
\begin{equation}
\psi_\nu(\theta) = 
-\frac{2}{\pi}i
\frac{\exp(i(\nu+1)\theta)}{\Gamma(\nu+1)}
\int_0^\infty  \int_{0}^\infty 
\frac{\exp(-(1+\beta x) w)}{\sqrt{x^2-2ix}}
 w^\nu
\ dx \ dw.
\label{integral:5}
\end{equation}
Since $\real(1+\beta x) >0$, we may use Formula~(\ref{preliminaries:laplace2})
to evaluate the integral with respect to $w$ in (\ref{integral:5})
whenever $\real(\nu) > -1$.    In this way we conclude that 
\begin{equation}
\psi_\nu(\theta)=
-\frac{2}{\pi}i
\exp(i(\nu+1)\theta)
\int_0^\infty  
\frac{1}{\sqrt{x^2-2ix}\left(1+\beta x\right)^{\nu+1}} 
\ dx,
\label{integral:psi}
\end{equation}
where $\beta = \sin(\theta)\exp(i\theta)$, 
for all $\real(\nu) > -1$ and $0 < \theta < \pi/2$.

\label{section:integral}
\end{section}

\begin{section}{A nonoscillatory phase function for Legendre's equation}

\begin{figure}[t!!]
\begin{center}
\includegraphics[width=.45\textwidth]{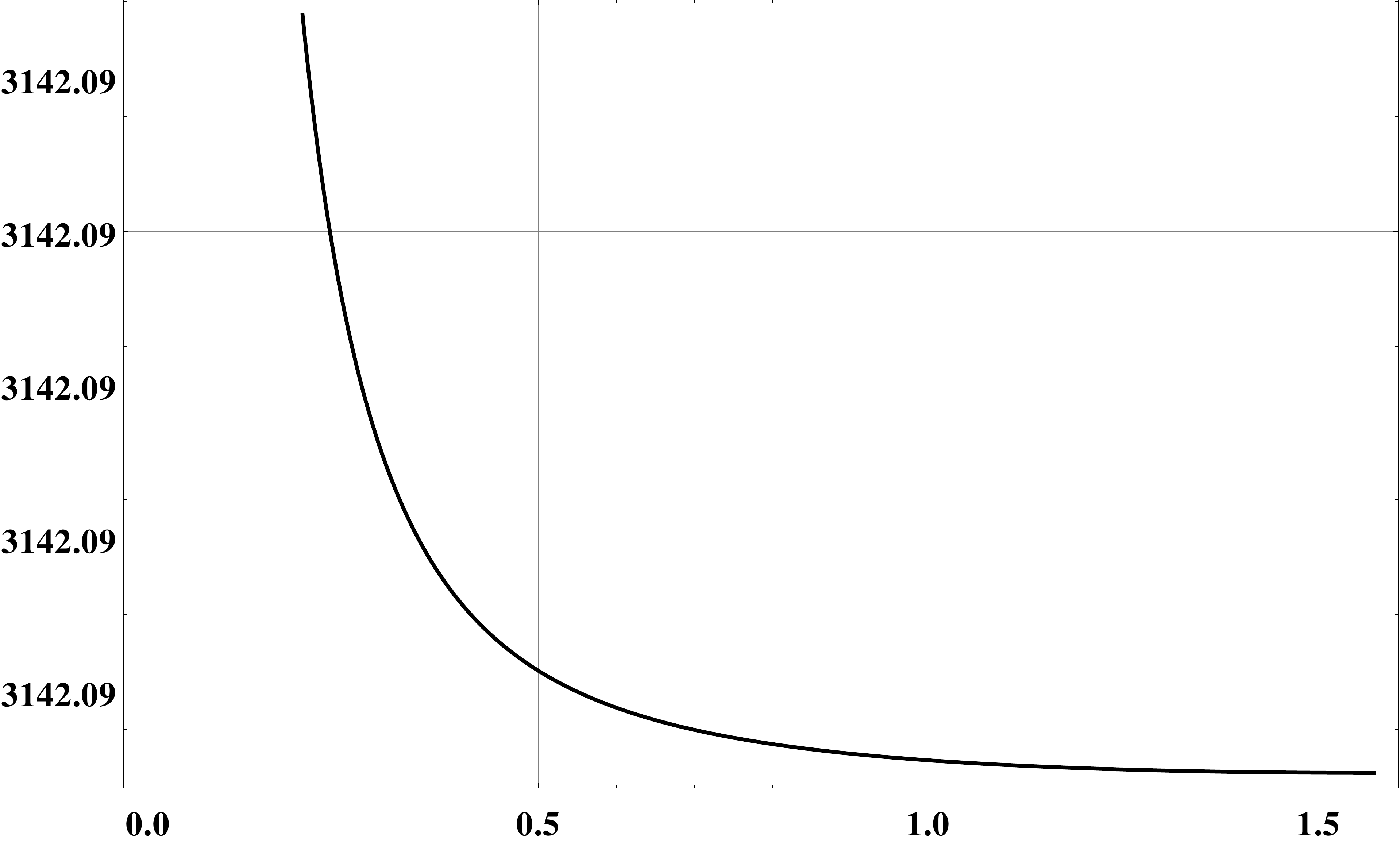}
\hfil
\includegraphics[width=.45\textwidth]{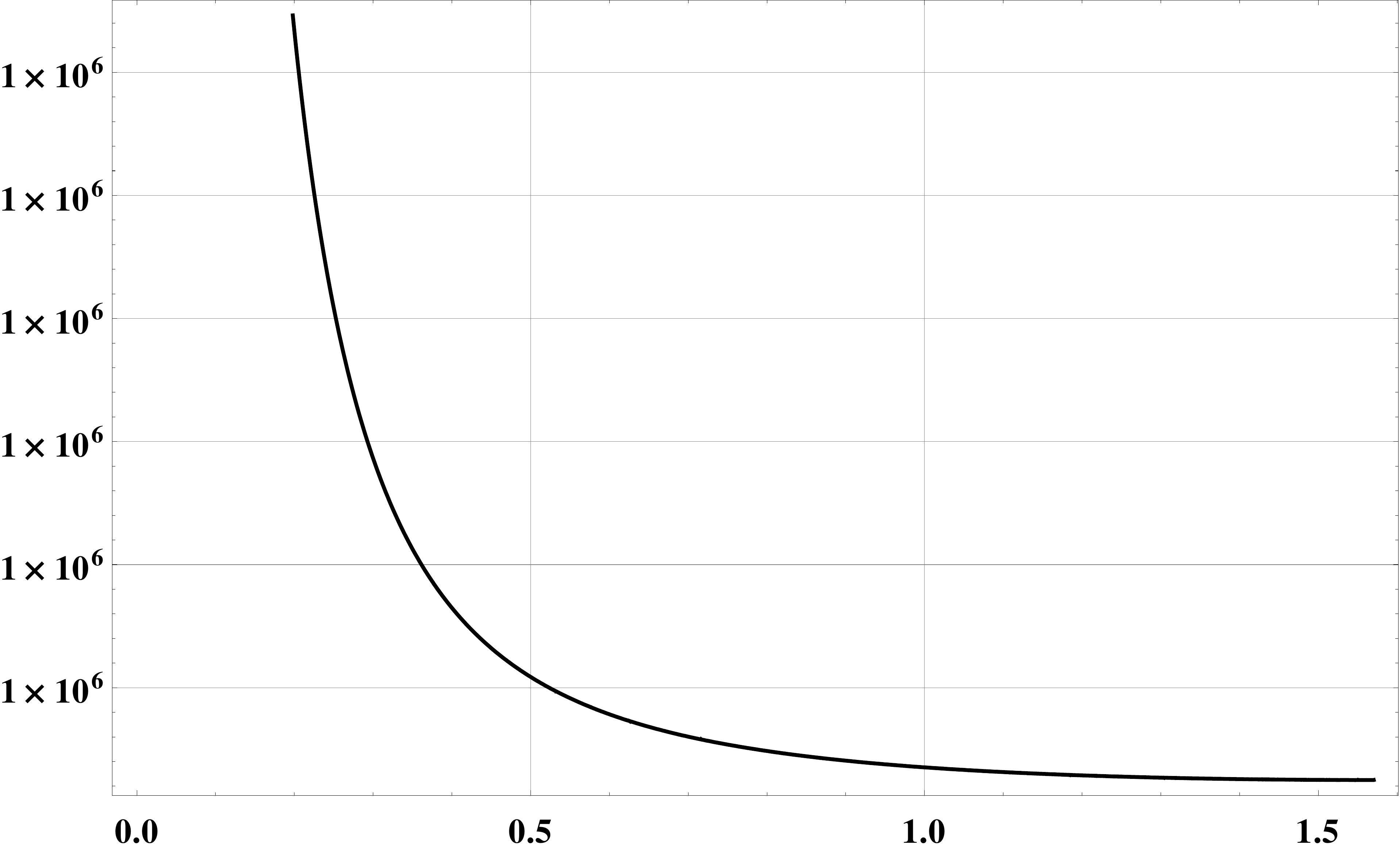}
\caption{
Plots of the imaginary part of the logarithmic derivative of $\psi_\nu$, which
is the derivative of a nonoscillatory phase function for Legendre's differential equation,
when $\nu = 100 \pi$ (left) and when $\nu=10^6$ (right).
}
\label{figure2}
\end{center}
\end{figure}

Formula~(\ref{integral:psi}) expresses $\psi_\nu$  as the product
of an oscillating exponential function and an  integral involving only
nonoscillatory elementary functions.
It follows that $\psi_\nu$ gives rise to a nonoscillatory phase
function for Legendre's differential equation.  To see this, 
we first define $\sigma_\nu$ via the formula
\begin{equation}
\sigma_\nu(\theta) = 
\int_0^\infty  
\frac{1}{\sqrt{x^2-2ix}\left(1+\sin(\theta) \exp(i\theta) x\right)^{\nu+1}} 
\ dx
\label{nonoscillatory:1}
\end{equation}
so that
\begin{equation}
\psi_\nu(\theta) = - \frac{2}{\pi}i
\exp( i (\nu+1) \theta) \sigma_\nu(\theta).
\label{nonoscillatory:2}
\end{equation}
Next, we observe that
\begin{equation}
\begin{aligned}
\sigma_\nu'(\theta) &= 
-(\nu+1)
\exp(2i\theta)
\int_0^\infty  
\frac{x}
{\sqrt{x^2-2ix}\left(1+\sin(\theta) \exp(i\theta) x\right)^{\nu+2}} 
\ dx.
\end{aligned}
\label{nonoscillatory:3}
\end{equation}
Now, we  combine (\ref{nonoscillatory:1}), (\ref{nonoscillatory:2}) and
(\ref{nonoscillatory:3}) to obtain
\begin{equation}
\begin{aligned}
\frac{\psi_\nu'(\theta)}{\psi_\nu(\theta)}
&=  
\frac
{ i (\nu+1) \exp(i(\nu+1)\theta)  \sigma_\nu(\theta)
+
\exp(i(\nu+1)\theta)  \sigma_\nu'(\theta)
}
{ \exp(i(\nu+1)\theta)\sigma_\nu(\theta)}
\\
&=
i (\nu+1) + 
\frac{  \sigma_\nu'(\theta)}
{ \sigma_\nu(\theta)}
\\
&=
i (\nu+1) 
- (\nu+1)
 \exp(2 i\theta)
\left(\int_0^\infty \frac{x\  dx}{\sqrt{x^2-2ix}\left(1+\sin(\theta) \exp(i\theta) x\right)^{\nu}} \right)
\cdot \\
&\ \ \ \ \ \ \ \ \ \ 
\left(\int_0^\infty \frac{dx}{\sqrt{x^2-2ix}\left(1+\sin(\theta) \exp(i\theta) x\right)^{\nu}} \right)^{-1}.
\end{aligned}
\label{nonoscillatory:alphap}
\end{equation}
According to Theorem~\ref{preliminaries:theorem1}, the imaginary part
of $\psi_\nu'/\psi_\nu$ 
 is the derivative of a phase function $\alpha_\nu$ for
Legendre's differential equation (\ref{preliminaries:legendre}).  We conclude form
Formula~(\ref{nonoscillatory:alphap})
 that $\alpha_\nu$ and $\alpha_\nu'$ are nonoscillatory functions.
Figure~\ref{figure2} contains plots of the imaginary part of the logarithmic
derivative of the function $\psi_\nu$ when $\nu = 100 \pi$ and when $\nu=10^6$.

\label{section:nonoscillatory}
\end{section}


\begin{section}{An asymptotic formula for Legendre functions of large degrees}
\label{section:asymptotics}

In this section, we derive an asymptotic expansion for the Legendre functions
of large degrees.  Our starting point is the formula
\begin{equation}
\psi_\nu(\theta)=
-\frac{2}{\pi} i \exp(i(\nu+1)\theta)
\int_0^\infty  
\frac{1}{\sqrt{\tau^2-2i \beta \tau}\left(1+\tau\right)^{\nu+1}} 
\ d\tau,
\label{asymptotic:1}
\end{equation}
which is obtained from (\ref{integral:psi}) by 
  introducing the new variable $\tau=\beta x = \sin(\theta) \exp(i\theta) x $.
We replace the function
\begin{equation}
f(\tau) = \frac{1}{(1+\tau)^{\nu+1}}
\end{equation}
in (\ref{asymptotic:1})  with a sum of the form
\begin{equation}
g(\tau) = a_0 \exp(- p \tau) + \sum_{k=1}^N a_k \exp( - (p + k q ) \tau) 
+ b_k \exp( - (p - k q ) \tau),
\label{asymptotic:2}
\end{equation}
where $p = \nu+1$, $q = \sqrt{p}$,  and the coefficients $a_0,a_1,b_1,a_2,b_2,\ldots,a_N,b_N$ are chosen so
 the power series expansions of $f$ and $g$ around $0$
agree to order $2N$. That is, we require that
the system of $2N+1$ linear equations 
\begin{equation}
f^{(k)}(0) = g^{(k)}(0) \ \ \mbox{for all}\ \ k=0,1,\ldots,2N
\label{asymptotic:3}
\end{equation}
in the $2N+1$ variables  $a_0,a_1,b_1,a_2,b_2,\ldots,a_N,b_N$
be satisfied.  By so doing  we obtain the approximation
\begin{equation}
\begin{aligned}
\psi_\nu(\theta) \approx
-\frac{2}{\pi} i \exp(i (\nu+1) \theta) \left( 
a_0 \int_0^\infty  
\right.&\left.\frac{\exp(-p\tau)}{\sqrt{\tau^2 - 2i\beta}}\ d\tau   +
\sum_{k=1}^N a_k
\int_0^\infty  
\frac{\exp(-\left(p+kq\right)\tau)}{\sqrt{\tau^2 - 2i\beta}}\ d\tau  
\right. \\
&+\left.
\sum_{k=1}^N b_k
\int_0^\infty  
\frac{\exp(-\left(p-kq\right)\tau)}{\sqrt{\tau^2 - 2i\beta}}\ d\tau  
\right).
\end{aligned}
\label{asymptotic:4}
\end{equation}
Now by applying Formula~(\ref{preliminaries:laplace1}) to (\ref{asymptotic:4}) 
we conclude that
\begin{equation}
\begin{aligned}
\psi_\nu&(\theta) \approx
  \exp(i (\nu+1) \theta) \bigg( 
a_0    \exp(- i \beta p) H_0(\beta p )  +
\\&     
\sum_{k=1}^N a_k
\exp(- i \beta \left( p + kq\right) ) H_0\left(\beta \left( p + kq\right) \right)
+
b_k \exp(- i \beta \left( p - kq\right) ) H_0\left(\beta \left( p - kq\right) \right)
\bigg),
\end{aligned}
\label{asymptotic:psi}
\end{equation}
where  $\beta = \sin(\theta) \exp(i\theta)$, $p = \nu+1$, and $q = \sqrt{p}$.
The use of Formula~(\ref{preliminaries:laplace1}) is justified so long
as $\real(p) > N^2$ and $0 < \theta < \pi/2$ (the second condition  ensures
that $|\arg(2i\beta)| < \pi$).   
The linear system (\ref{asymptotic:3})  can be solved easily using a
 computer algebra system.  Our Mathematica script for doing so
appears in an appendix of  this paper.  A second appendix lists
the coefficients when $N=2$, $N=3$, $N=4$, $N=5$  and $N=6$.

Error bounds for the formula (\ref{asymptotic:psi})
can be derived quite easily using standard techniques (the monograph \cite{Olver}, for example,
contains many similar estimates).  We omit them here, though, 
because  our point is not that (\ref{asymptotic:psi}) is especially accurate,
 but rather that its form is quite different  from most widely used
asymptotic expansions for Legendre functions in that it expresses $\psi_\nu$ as the product
of an  oscillatory exponential function and a sum of nonoscillatory
functions rather than as a sum of oscillatory functions.   
We do, however, report on extensive numerical experiments
which assess the performance of 
Formula~(\ref{asymptotic:psi}) in Section~\ref{section:numerics} of this article.
That the functions
\begin{equation}
\exp\left(- i \beta \left( p \pm kq\right) \right) H_0\left(\beta \left( p \pm kq\right)\right)
\end{equation}
appearing  in (\ref{asymptotic:psi}) are nonoscillatory is obvious from the formula
\begin{equation}
\exp\left(-iz\right)
H_0\left(z \right)
= -\frac{2}{\pi}i\int_0^\infty \frac{\exp(-z  x)}{\sqrt{x^2-2ix}}\ dx 
\end{equation}
obtained by letting $b=-2i$ in (\ref{preliminaries:laplace1}).
Stieltjes' formula (\ref{preliminaries:stieltjes}),
which approximates $P_\nu$ via a sum of cosines which oscillate at frequencies
on the order of $\nu$, furnishes an example of the latter type of expansion,
as does  Olver's asymptotic expansion
\begin{equation}
P_n(\cos(\theta)) \approx
\left(\frac{\theta}{\sin(\theta)}\right)^{1/2}
\left(
J_0(u\theta) \sum_{k=0}^m \frac{A_k(\theta)}{u^{2k}} +
\frac{\theta}{u} J_{-1}(u \theta) 
\sum_{k=0}^m \frac{B_k(\theta)}{u^{2k}}
\right),
\label{asymptotic:Olver}
\end{equation}
where $u= n +1/2$.   The definition of the coefficients
$A_k$ and $B_k$ as well as a derivation of (\ref{asymptotic:Olver}) can be
found in Chapter~12 of \cite{Olver}.  Yet another example
is given by Dunster's convergent expansions
 for Legendre functions \cite{Dunster-Legendre2}, which have a form
similar to (\ref{asymptotic:Olver}).

Among other things, the form of the expansion (\ref{asymptotic:psi})
is conducive to computing the derivative of the 
nonoscillatory phase function $\alpha_\nu$
associated with $\psi_\nu$.  To see this, we rewrite (\ref{asymptotic:psi}) as 
\begin{equation}
\begin{aligned}
\psi_\nu(\theta)  \approx
 \exp(i (\nu+1) \theta) \bigg( 
a_0    S (\beta p) 
  +      \sum_{k=1}^N a_k S( (p+kq)\beta) 
+
b_k S(\beta \left( p - kq\right) ) \bigg),
\end{aligned}
\end{equation}
where the function $S$ is defined by
\begin{equation}
S(z) = \exp(-iz) H_0(z),
\end{equation}
and differentiate both sides in order to obtain
\begin{equation}
\begin{aligned}
&\psi_\nu'(\theta)  \approx
i(\nu+1) \exp(i (\nu+1) \theta) \bigg(  a_0   S (\beta p)   +      \sum_{k=1}^N a_k S( (p+kq)\beta) +
b_k S(\beta \left( p - kq\right) ) \bigg)  \\
&+
\exp(i (\nu+3) \theta) \bigg(  a_0  p  S' (\beta p) \sum_{k=1}^N a_k (p+kq) S'( (p+kq)\beta) 
 + b_k (p-kq) S'(\beta \left( p - kq\right) ) \bigg).  \\
\end{aligned}
\label{asymptotic:psider}
\end{equation}
From (\ref{preliminaries:h0der}) we see that
\begin{equation}
S'(z) = \frac{d}{dz} \left( \exp(-iz) H_0(z) \right) = 
-\exp(-iz) H_1(z) - i\exp(-iz)H_0(z).
\end{equation}
By combining (\ref{asymptotic:psider}) and (\ref{asymptotic:psi}) we can evaluate
the logarithmic derivative $\psi_\nu'/ \psi_\nu$ of $\psi_\nu$ and hence
the derivative of a nonoscillatory phase function for Legendre's differential equation,
which is the imaginary part of this ratio.

\label{section:asymptotic}
\end{section}

\begin{section}{Numerical experiments}
\label{section:numerics}

In this section, we describe numerical experiments which were conducted to assess
the performance of the asymptotic expansions of Section~\ref{section:asymptotics}.
Our code was written in Fortran and compiled with the GNU Fortran Compiler
version 5.2.1.  The calculations were carried out on a laptop
equipped with an Intel i7-5600U processor running at  2.60GHz.

\begin{subsection}{The accuracy of the expansion (\ref{asymptotic:psi}) as a function of $N$ and $\nu$}

We measured the accuracy of the expansion (\ref{asymptotic:psi}) for various values of $N$
and $\nu$.
For each pair of values of $N$ and $\nu$ considered, we evaluated (\ref{asymptotic:psi}) at a collection
of $1\sep{,}000$ points  on the interval $(0,\pi/2)$.  The first $500$ points 
were drawn at random from the uniform distribution on the interval $(0,\pi/2)$, while the remaining
 points were constructed by drawing $500$ points from the uniform
distribution on the interval $(0,1)$ and applying the mapping
$t \to \exp(-36 t)$
to each of them.  In this way, we ensured that the accuracy of (\ref{asymptotic:psi}) was tested 
near the singular point of Legendre's equation which occurs when $\theta = 0$.

The results are reported in Tables~\ref{table1} and \ref{table2}.  
 Table~\ref{table1} gives the maximum relative error in the value of $\psi_\nu$ 
which was observed when these calculations were carried out in double precision
arithmetic as a function of $\nu$ and $N$, and 
Table~\ref{table2} reports the maximum relative error in the value of $\psi_\nu$ 
which was observed when they were performed using quadruple precision arithmetic
as a function of $\nu$ and $N$.  Reference values
were computed by running the algorithm of \cite{BremerKummer} in quadruple precision
arithmetic.    Like the asymptotic expansion (\ref{asymptotic:psi}), the
algorithm of \cite{BremerKummer} is capable of achieving high accuracy
near the singular point of Legendre's equation
whereas
 Stieltjes' expansion (\ref{preliminaries:stieltjes}) and the well-known
three term recurrence relations are inaccurate when $\theta$ is near 0.
See \cite{Bogaert-Michiels-Fostier, Hale-Townsend}, though, for
the derivation of asymptotic expansions of Legendre functions
which are accurate for $\theta$ near 0 and \cite{Olver,Dunster-Legendre1,Dunster-Legendre2}
for expansions of Legendre functions which are accurate for all $\theta \in (0,\pi/2)$.

The routine we used to evaluate the Hankel function
of order $0$ achieves roughly double precision accuracy, even when
executed using quadruple precision arithmetic.  Thus the minimum
error achieved when the computations were performed using quadruple precision arithmetic
was on the order of $10^{-16}$ (see Table~\ref{table2}).

\begin{table}[t!!!]
\begin{center}
\begin{tabular}{c@{\hspace{3em}}ccccc}
\toprule
\addlinespace[.5em]
$\nu$ & $N=2$ & $N=3$ & $N=4$ & $N=5$ & $N=6$ \\
\midrule
\addlinespace[.75em]
$10^2$ & $1.55\e{-06}$ & $5.30\e{-08}$ & $1.48\e{-09}$ & $6.05\e{-11}$ & $1.17\e{-11}$ \\ 
\addlinespace[.25em]
$10^2 \pi $ & $5.02\e{-08}$ & $5.00\e{-10}$ & $2.84\e{-12}$ & $6.49\e{-14}$ & $5.63\e{-14}$ \\ 
\addlinespace[.25em]
$10^3 $ & $1.55\e{-09}$ & $4.74\e{-12}$ & $2.09\e{-13}$ & $2.09\e{-13}$ & $2.09\e{-13}$ \\ 
\addlinespace[.25em]
$10^3 \pi $ & $5.02\e{-11}$ & $1.16\e{-12}$ & $1.16\e{-12}$ & $1.16\e{-12}$ & $1.16\e{-12}$ \\ 
\addlinespace[.25em]
$10^4$ & $2.46\e{-12}$ & $1.90\e{-12}$ & $1.90\e{-12}$ & $1.90\e{-12}$ & $1.90\e{-12}$ \\ 
\addlinespace[.25em]
$10^4\pi$ & $6.70\e{-12}$ & $6.70\e{-12}$ & $6.70\e{-12}$ & $6.70\e{-12}$ & $6.70\e{-12}$ \\ 
\addlinespace[.25em]
$10^5$ & $2.17\e{-11}$ & $2.17\e{-11}$ & $2.17\e{-11}$ & $2.17\e{-11}$ & $2.17\e{-11}$ \\ 
\addlinespace[.25em]
$10^6$ & $2.15\e{-10}$ & $2.15\e{-10}$ & $2.15\e{-10}$ & $2.15\e{-10}$ & $2.15\e{-10}$ \\ 
\addlinespace[.25em]
$10^7$ & $2.00\e{-09}$ & $2.00\e{-09}$ & $2.00\e{-09}$ & $2.00\e{-09}$ & $2.00\e{-09}$ \\ 
\addlinespace[.25em]
$10^8$ & $2.33\e{-08}$ & $2.33\e{-08}$ & $2.33\e{-08}$ & $2.33\e{-08}$ & $2.33\e{-08}$ \\ 
\addlinespace[.25em]
$10^9$ & $2.15\e{-07}$ & $2.15\e{-07}$ & $2.15\e{-07}$ & $2.15\e{-07}$ & $2.15\e{-07}$ \\ 
\addlinespace[.25em]
\bottomrule
\end{tabular}
\end{center}
\caption{
The relative accuracy of the expansion (\ref{asymptotic:psi}) as a function 
 of $\nu$ and $N$. Here, the expansion was evaluated
using double precision arithmetic.}
\label{table1}
\end{table}

\begin{table}[b!!!!]
\begin{center}
\begin{tabular}{c@{\hspace{3em}}ccccc}
\toprule
\addlinespace[.5em]
$\nu$ & $N=2$ & $N=3$ & $N=4$ & $N=5$ & $N=6$ \\
\midrule
\addlinespace[.75em]
$10^2$ & $1.55\e{-06}$ & $5.30\e{-08}$ & $1.48\e{-09}$ & $6.05\e{-11}$ & $1.17\e{-11}$ \\ 
\addlinespace[.25em]
$10^2 \pi $ & $5.02\e{-08}$ & $5.00\e{-10}$ & $2.84\e{-12}$ & $3.28\e{-14}$ & $3.15\e{-15}$ \\ 
\addlinespace[.25em]
$10^3 $ & $1.55\e{-09}$ & $4.74\e{-12}$ & $7.06\e{-15}$ & $3.94\e{-16}$ & $4.10\e{-16}$ \\ 
\addlinespace[.25em]
$10^3 \pi $ & $5.02\e{-11}$ & $4.83\e{-14}$ & $3.72\e{-16}$ & $3.58\e{-16}$ & $3.61\e{-16}$ \\ 
\addlinespace[.25em]
$10^4$ & $1.55\e{-12}$ & $4.80\e{-16}$ & $4.94\e{-16}$ & $4.70\e{-16}$ & $4.73\e{-16}$ \\ 
\addlinespace[.25em]
$10^4\pi$ & $5.02\e{-14}$ & $4.95\e{-16}$ & $5.06\e{-16}$ & $4.84\e{-16}$ & $4.87\e{-16}$ \\ 
\addlinespace[.25em]
$10^5$ & $1.53\e{-15}$ & $3.41\e{-16}$ & $3.50\e{-16}$ & $3.28\e{-16}$ & $3.31\e{-16}$ \\ 
\addlinespace[.25em]
$10^6$ & $5.10\e{-16}$ & $4.86\e{-16}$ & $4.95\e{-16}$ & $4.73\e{-16}$ & $4.76\e{-16}$ \\ 
\addlinespace[.25em]
$10^7$ & $4.85\e{-16}$ & $4.57\e{-16}$ & $4.67\e{-16}$ & $4.42\e{-16}$ & $4.45\e{-16}$ \\ 
\addlinespace[.25em]
$10^8$ & $4.43\e{-16}$ & $4.14\e{-16}$ & $4.24\e{-16}$ & $3.98\e{-16}$ & $4.02\e{-16}$ \\ 
\addlinespace[.25em]
$10^9$ & $4.91\e{-16}$ & $4.63\e{-16}$ & $4.72\e{-16}$ & $4.48\e{-16}$ & $4.51\e{-16}$ \\ 
\addlinespace[.25em]
\bottomrule
\end{tabular}
\end{center}
\caption{
The relative accuracy of the expansion (\ref{asymptotic:psi}) as a function
 of $\nu$ and $N$.  Here, the expansion was evaluated
using quadruple precision arithmetic.}
\label{table2}
\end{table}

We also observe that relative accuracy
was lost as $\nu$ increases when (\ref{asymptotic:psi}) is  evaluated using  double
precision arithmetic.    This loss of precision is unsurprising and consistent
with the condition number of the evaluation of the highly oscillatory functions
$P_\nu$ and $Q_\nu$.
Indeed, evaluating $P_\nu(x)$ and $Q_\nu(x)$  is analogous to calculating
the values of $\cos(\nu x)$ and $\sin(\nu x)$.
If $\cos(\nu x)$ and $\sin(\nu x)$ can be evaluated with accuracy on the order of $\epsilon$,
then so can the value of $\mbox{mod}(\nu,1)$.  This clearly limits
the precision with which $\cos(\nu x)$ and $\sin(\nu x)$ can be evaluated using
finite precision arithmetic.  A similar argument applies to $P_\nu(x)$
and $Q_\nu(x)$.
The only exception  is when the value of $\mbox{mod}(\nu,1)$ is known to high precision 
(for instance, when $\nu$ is an integer).  Then $\cos(\nu x)$ and $\sin (\nu x)$
can be calculated with comparable precision, 
as can $P_\nu$ and $Q_\nu$  (via the well-known three term recurrence relations, for instance).

\label{section:numerics:accuracy}
\end{subsection}

\begin{subsection}{The accuracy of the expansion (\ref{asymptotic:psi}) as a function of $\theta$}

In order to measure the accuracy of the expansion (\ref{asymptotic:psi})
as a function of $\theta$, we sampled a collection of $1\sep{,}000$ points
on the interval $(0,\pi/2)$ using the same method as in 
Section~\ref{section:numerics:accuracy} and then evaluated (\ref{asymptotic:psi})
at each of the chosen values of $\theta$ 
with $\nu=1\sep{,}000$ and $N=2$, $N=3$ and $N=4$.  The base-$10$ logarithm
of the relative error in each of the resulting values was computed and a plot showing these errors
appears in Figure~\ref{figure1}.  Once again, reference values were computed by running
the algorithm of \cite{BremerKummer} in quadruple precision arithmetic.
We observe that the error in Formula~(\ref{asymptotic:psi}) is remarkably uniform as a function
of $\theta$ --- it is  nearly constant until $\theta$ nears the singular point at $0$,
 at which point it decreases slightly.

\begin{figure}[h!!!]
\begin{center}
\includegraphics[width=.8\textwidth]{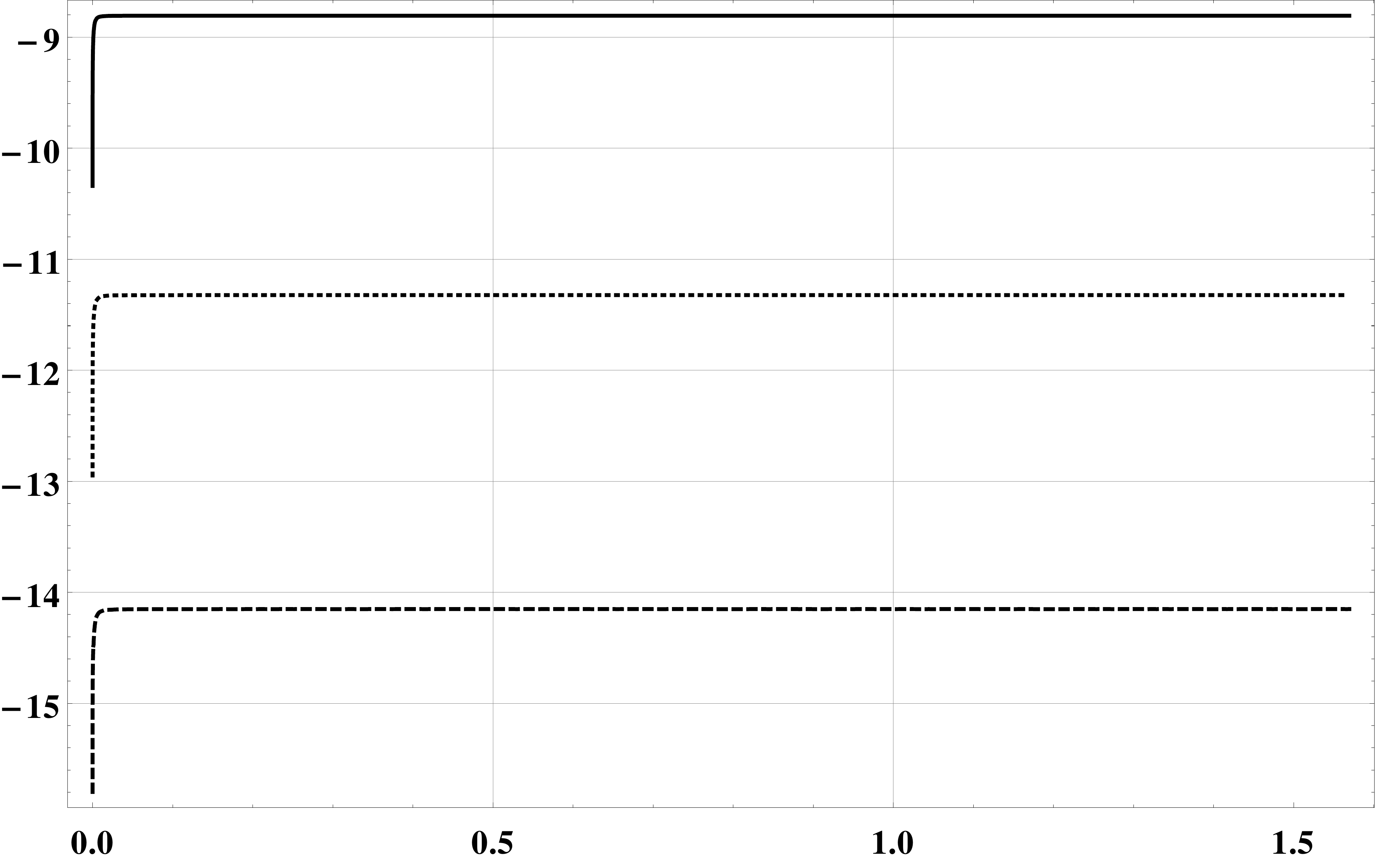}
\caption{\small
The base-$10$ logarithm of the
relative accuracy of the expansion
(\ref{asymptotic:psi}) as a function of $\theta$ when $\nu = 1\sep{,}000$
and $N=2$ (top line), $N=3$ (middle line), $N=4$ (bottom line).
}
\label{figure1}
\end{center}
\end{figure}

\end{subsection}

\begin{subsection}{The speed of the expansion (\ref{asymptotic:psi}) as function of $N$}

Next, we  measured the time required to
evaluate $\psi_\nu$ using the expansion (\ref{asymptotic:psi}).
In particular, for several pairs of values  of $N$ and $\nu$, we evaluated
(\ref{asymptotic:psi}) at a collection of $1\sep{,}000$ points drawn
from the uniform distribution on the interval $(0,\pi/2)$.
We also applied the same procedure to Stieltjes' expansion (\ref{preliminaries:stieltjes})
with $M=16$.   Table~\ref{table:3} gives the average time required
to evaluate (\ref{asymptotic:psi}) and (\ref{preliminaries:stieltjes}).

We note that while the time required to evaluate (\ref{asymptotic:psi})
is slightly larger than the time required to evaluate
Stieltjes' formula (\ref{preliminaries:stieltjes}),  the expansion
of Section~\ref{section:asymptotics} gives both the value of $P_\nu$ 
and that of $Q_\nu$ while Stieltjes' formula gives only the value of $P_\nu$.  
Moreover, unlike
Stieltjes' formula, our asymptotic approximation is accurate for $\theta$ near $0$.

\begin{table}[h!]
\begin{center}
\begin{tabular}{c@{\hspace{3em}}ccccc}
\toprule
\addlinespace[.5em]
$\nu$ & $N=3$ & $N=4$ & $N=5$ & $N=6$ & Stieltjes' formula \\
\midrule
\addlinespace[.75em]
$10^2$ & $2.36\e{-06}$ & $2.93\e{-06}$ & $3.92\e{-06}$ & $4.27\e{-06}$ & $1.57\e{-06}$ \\ 
\addlinespace[.25em]
$10^2 \pi$ & $2.35\e{-06}$ & $3.02\e{-06}$ & $3.98\e{-06}$ & $4.31\e{-06}$ & $1.71\e{-06}$ \\ 
\addlinespace[.25em]
$10^3$ & $2.38\e{-06}$ & $2.98\e{-06}$ & $3.95\e{-06}$ & $4.31\e{-06}$ & $1.58\e{-06}$ \\ 
\addlinespace[.25em]
$10^3 \pi$ & $2.38\e{-06}$ & $2.96\e{-06}$ & $3.93\e{-06}$ & $4.26\e{-06}$ & $1.64\e{-06}$ \\ 
\addlinespace[.25em]
$10^4$ & $2.38\e{-06}$ & $3.00\e{-06}$ & $3.94\e{-06}$ & $4.41\e{-06}$ & $1.55\e{-06}$ \\ 
\addlinespace[.25em]
$10^4 \pi$ & $2.33\e{-06}$ & $2.97\e{-06}$ & $3.95\e{-06}$ & $4.26\e{-06}$ & $1.56\e{-06}$ \\ 
\addlinespace[.25em]
$10^5$ & $2.43\e{-06}$ & $2.95\e{-06}$ & $3.93\e{-06}$ & $4.31\e{-06}$ & $1.63\e{-06}$ \\ 
\addlinespace[.25em]
$10^6$ & $2.33\e{-06}$ & $2.94\e{-06}$ & $3.93\e{-06}$ & $4.23\e{-06}$ & $1.55\e{-06}$ \\ 
\addlinespace[.25em]
$10^7$ & $2.44\e{-06}$ & $2.98\e{-06}$ & $3.96\e{-06}$ & $4.37\e{-06}$ & $1.54\e{-06}$ \\ 
\addlinespace[.25em]
$10^8$ & $2.36\e{-06}$ & $3.01\e{-06}$ & $3.99\e{-06}$ & $4.23\e{-06}$ & $1.61\e{-06}$ \\ 
\addlinespace[.25em]
$10^9$ & $2.39\e{-06}$ & $2.99\e{-06}$ & $3.94\e{-06}$ & $4.38\e{-06}$ & $1.67\e{-06}$ \\ 
\addlinespace[.25em]
\bottomrule
\end{tabular}
\end{center}
\caption{
A comparison of the average time (in seconds) required to evaluate (\ref{asymptotic:psi})
for various values of $N$ and $\nu$ with the time requires to evaluate the first $17$ terms
of Stieltjes' expansion (\ref{preliminaries:stieltjes}).  Note that
Stieltjes' formula yields only the value of $P_\nu$ while (\ref{asymptotic:psi})
yields both $P_\nu$ and $Q_\nu$.  Moreover, unlike (\ref{asymptotic:psi}),
Stieltjes' is not accurate for arguments close to the singular points of Legendre's differential
equation.
}
\label{table:3}
\end{table}

\end{subsection}

\begin{subsection}{Evaluation of the derivative of a phase function for Legendre's equation}

In this last experiment, we measured the accuracy achieved
when Formulas~(\ref{asymptotic:psi}) and (\ref{asymptotic:psider})
are combined in order to evaluate the derivative of the nonoscillatory phase function $\alpha_\nu$
for Legendre's equation associated with the function $\psi_\nu$.
More specifically, for each of 
several pairs of values of $N$ and $\nu$, we evaluated
(\ref{asymptotic:psi}) and (\ref{asymptotic:psider}) at a collection of $1\sep{,}000$
points in the interval $(0,\pi/2)$ and used these quantities to calculate
$\alpha_\nu'$ at each of these points.  The points were chosen as in 
Section~\ref{section:numerics:accuracy}.   These experiments were performed using
double precision arithmetic.   The obtained values of $\alpha_\nu'$ were compared
with reference values computed by running the algorithm of \cite{BremerKummer}
using quadruple precision arithmetic.

The results are reported in Table~\ref{table4}.   For each pair of values
of $\nu$ and $N$, it lists the maximum relative error in $\alpha_\nu'$ which
was observed.  We note that, unlike the experiments of Section~\ref{section:numerics:accuracy},
near double precision accuracy was obtained by performing the calculations
in double precision arithmetic.  This is not surprising since the condition
number of the evaluation of the nonoscillatory function
$\alpha_\nu'$ is small and not dependent on $\nu$.

\begin{table}[h!]
\begin{center}
\begin{tabular}{c@{\hspace{3em}}ccccc}
\toprule
\addlinespace[.5em]
$\nu$ & $N=2$ & $N=3$ & $N=4$ & $N=5$ & $N=6$ \\
\midrule
\addlinespace[.75em]
$10^2$ & $4.87\e{-07}$ & $2.07\e{-08}$ & $7.36\e{-09}$ & $7.04\e{-09}$ & $7.03\e{-09}$ \\ 
\addlinespace[.25em]
$10^2 \pi $ & $1.53\e{-08}$ & $1.64\e{-10}$ & $1.60\e{-10}$ & $1.60\e{-10}$ & $1.60\e{-10}$ \\ 
\addlinespace[.25em]
$10^3$ & $4.78\e{-10}$ & $1.22\e{-12}$ & $2.42\e{-15}$ & $1.27\e{-15}$ & $1.41\e{-15}$ \\ 
\addlinespace[.25em]
$10^3 \pi$ & $1.54\e{-11}$ & $1.29\e{-14}$ & $1.07\e{-15}$ & $1.27\e{-15}$ & $1.36\e{-15}$ \\ 
\addlinespace[.25em]
$10^4 $ & $4.78\e{-13}$ & $1.36\e{-15}$ & $1.36\e{-15}$ & $1.48\e{-15}$ & $1.36\e{-15}$ \\ 
\addlinespace[.25em]
$10^4 \pi $ & $1.56\e{-14}$ & $1.08\e{-15}$ & $1.26\e{-15}$ & $1.08\e{-15}$ & $1.34\e{-15}$ \\ 
\addlinespace[.25em]
$10^5 $ & $1.10\e{-15}$ & $9.65\e{-16}$ & $1.38\e{-15}$ & $1.14\e{-15}$ & $1.20\e{-15}$ \\ 
\addlinespace[.25em]
$10^6 $ & $1.30\e{-15}$ & $1.30\e{-15}$ & $1.19\e{-15}$ & $1.30\e{-15}$ & $1.33\e{-15}$ \\ 
\addlinespace[.25em]
$10^7 $ & $1.22\e{-15}$ & $1.22\e{-15}$ & $1.22\e{-15}$ & $1.70\e{-15}$ & $1.62\e{-15}$ \\ 
\addlinespace[.25em]
$10^8 $ & $1.56\e{-15}$ & $1.42\e{-15}$ & $1.45\e{-15}$ & $1.56\e{-15}$ & $1.42\e{-15}$ \\ 
\addlinespace[.25em]
$10^9 $ & $1.12\e{-15}$ & $1.34\e{-15}$ & $1.34\e{-15}$ & $1.41\e{-15}$ & $1.38\e{-15}$ \\ 
\addlinespace[.25em]
\bottomrule
\end{tabular}
\end{center}
\caption{
The relative accuracy with which the derivative of the nonoscillatory phase function
for Legendre's differential equation is evaluated via
the expansions of Section~\ref{section:asymptotics}.
These calculations were performed using double precision arithmetic.
}
\label{table4}
\end{table}


\end{subsection}

\end{section}

\begin{section}{Conclusions}
\label{section:conclusion}

Nonoscillatory phase functions are powerful analytic and numerical tools.
Among other things, explicit formulas for them can be used to efficiently and accurately
evaluate special functions, their zeros and to apply special function transforms.

Here, we factored a particular solution of Legendre's differential equation
as the product of an oscillatory exponential function and an integral
involving only nonoscillatory elementary functions.
 By so doing, we proved the existence of a nonoscillatory phase
function and derived an asymptotic formula  of 
an unusual type.  Specifically, our formula represents the oscillatory
functions $P_\nu$ and $Q_\nu$
as the product of an oscillatory exponential function 
 and a sum of nonoscillatory functions.

We will report on the use of the results of the paper to apply 
the Legendre transform rapidly and on
 generalizations of this work to the case of associated
Legendre functions, prolate spheroidal wave functions and other related
special functions at a later date.

\end{section}

\begin{section}{Acknowledgments}
James Bremer was supported  by National Science Foundation grant DMS-1418723.
Vladimir Rokhlin was supported by the Office of Naval Research under contract N00014-16-1-2123
and by the Air Force Scientific Research Office under contract FA9550-16-1-0175.
\end{section}


\vfill\eject


\appendix

\begin{section}{Mathematica source code}

Below is a Mathematica script for generating the coefficients 
$a_0,a_1,b_1,\ldots,a_n,b_n$  in the  expansion (\ref{asymptotic:psi}).
\vskip 1em

{
\small
\begin{verbatim}
(***************************************************************************************)
ClearAll["Global`*"]

(* Construct the (2n+1)-term expansion *)
n=3;
f[t_]=1/(1+t)^(q^2);
g[t_]=a0*Exp[-q^2*t]+Sum[a[k]*Exp[-(q^2+k*q)*t]+b[k]*Exp[-(q^2-k*q)*t],{k,1,n}];

eqs={};
vars={};
h0[t_]=f[t]-g[t];
For[k=0,k<=2*n,k++,eqs=Simplify[Join[eqs,{h0[0]==0}]];h0[t_]=h0'[t]];
For[k=1,k<=n,k++,vars=Join[vars,{a[k],b[k]}]];
vars=Join[vars,{a0}];
sol=Simplify[Solve[eqs,vars][[1]],p>0]

(* Write the TEX version of the expansion to the file tmp2 *)
OpenWrite["tmp2"];
WriteString["tmp2","\\begin{equation}\n\\begin{aligned}\n"];
WriteString["tmp2","a_0 &="ToString[TeXForm[a0/.sol]],", \\\\ \n"];
For[k=1,k<=n,k++,
WriteString["tmp2","a_",k," &="ToString[TeXForm[a[k]/.sol]],", \\\\ \n"]];
For[k=1,k<=n-1,k++,
WriteString["tmp2","b_",k," &="ToString[TeXForm[b[k]/.sol]],", \\\\ \n"]];
WriteString["tmp2","b_",k," &="ToString[TeXForm[b[k]/.sol]],".\n"];
WriteString["tmp2","\\end{aligned}\n\\end{equation}"];
Close["tmp2"];

(* Write the FORTRAN version of the expansion to tmp *)
OpenWrite["tmp"];
WriteString["tmp","a0 = ",ToString[FortranForm[N[Expand[a0/.sol],36]]],"\n"];
For[k=1,k<=n,k++,
WriteString["tmp","a(",k,") ="ToString[FortranForm[N[Expand[a[k]/.sol],36]]],"\n"]];
For[k=1,k<=n-1,k++,
WriteString["tmp","b(",k,") ="ToString[FortranForm[N[Expand[b[k]/.sol],36]]],"\n"]];
WriteString["tmp","b(",k,") ="ToString[FortranForm[N[Expand[b[k]/.sol],36]]],"\n"];
Close["tmp"];
(***************************************************************************************)
\end{verbatim}
}

\end{section}


\vfill\eject

\begin{section}{The coefficients in the expansion (\ref{asymptotic:psi})}

When $N=2$, the coefficients are
\begin{equation*}
\tiny
\begin{aligned}
a_0 &= \frac{3}{2 q^2}+\frac{1}{2}, \\ 
a_1 &= \frac{q^2-2 q-6}{6 q^2}, \\ 
a_2 &= \frac{q^2+2 q+3}{12 q^2}, \\ 
b_1 &= \frac{q^2+2 q-6}{6 q^2}, \\ 
b_2 &= \frac{q^2-2 q+3}{12 q^2}.
\end{aligned}
\end{equation*}
%

When $N=3$, they are
\begin{equation*}
\tiny
\begin{aligned}
a_0 &= -\frac{-7 q^4+23 q^2+60}{18 q^4}, \\ 
a_1 &= \frac{6 q^4-3 q^3+26 q^2+12 q+60}{24 q^4}, \\ 
a_2 &= -\frac{-3 q^4+35 q^2+24 q+60}{60 q^4}, \\ 
a_3 &= \frac{2 q^4+15 q^3+50 q^2+36 q+60}{360 q^4}, \\ 
b_1 &= \frac{6 q^4+3 q^3+26 q^2-12 q+60}{24 q^4}, \\ 
b_2 &= \frac{3 q^4-35 q^2+24 q-60}{60 q^4}, \\ 
b_3 &= \frac{2 q^4-15 q^3+50 q^2-36 q+60}{360 q^4}.
\end{aligned}
\end{equation*}

When $N=4$:
\begin{equation*}
\tiny
\begin{aligned}
\\
a_0 &=\frac{115 q^6+59 q^4+1854 q^2+2520}{288 q^6}, \\
a_1 &=-\frac{-87 q^6+59 q^5+37 q^4+114 q^3+1914 q^2+360 q+2520}{360 q^6}, \\ 
a_2 &=\frac{39 q^6+28 q^5+7 q^4+300 q^3+2094 q^2+720 q+2520}{720 q^6}, \\ 
a_3 &=-\frac{-11 q^6-63 q^5+77 q^4+630 q^3+2394 q^2+1080 q+2520}{2520 q^6}, \\ 
a_4 &=\frac{3 q^6+56 q^5+427 q^4+1176 q^3+2814 q^2+1440 q+2520}{20160 q^6}, \\ 
b_1 &=\frac{87 q^6+59 q^5-37 q^4+114 q^3-1914 q^2+360 q-2520}{360 q^6}, \\ 
b_2 &=\frac{39 q^6-28 q^5+7 q^4-300 q^3+2094 q^2-720 q+2520}{720 q^6}, \\ 
b_3 &=-\frac{-11 q^6+63 q^5+77 q^4-630 q^3+2394 q^2-1080 q+2520}{2520 q^6}, \\ 
b_4 &=\frac{3 q^6-56 q^5+427 q^4-1176 q^3+2814 q^2-1440 q+2520}{20160 q^6}.
\end{aligned}
\end{equation*}

When $N=5$:
\begin{equation*}
\tiny
\begin{aligned}
\\
a_0 &=-\frac{-359 q^8+20 q^6+1530 q^4+21636 q^2+22680}{900 q^8}, \\ 
a_1 &=\frac{2091 q^8-1396 q^7+747 q^6-104 q^5+12654 q^4+12672 q^3+175608 q^2+20160 q+181440}{8640 q^8}, \\ 
a_2 &=-\frac{-204 q^8-137 q^7+372 q^6-259 q^5+3654 q^4+6876 q^3+45792 q^2+10080 q+45360}{3780 q^8}, \\ 
a_3 &=\frac{179 q^8+1068 q^7+403 q^6-2184 q^5+20286 q^4+46656 q^3+195768 q^2+60480 q+181440}{40320 q^8}, \\ 
a_4 &=-\frac{-3 q^8-53 q^7-276 q^6-7 q^5+4158 q^4+9036 q^3+26676 q^2+10080 q+22680}{22680 q^8}, \\ 
a_5 &=\frac{3 q^8+100 q^7+1635 q^6+13160 q^5+58590 q^4+106560 q^3+236088 q^2+100800 q+181440}{1814400 q^8}, \\ 
b_1 &=\frac{2091 q^8+1396 q^7+747 q^6+104 q^5+12654 q^4-12672 q^3+175608 q^2-20160 q+181440}{8640 q^8}, \\ 
b_2 &=-\frac{-204 q^8+137 q^7+372 q^6+259 q^5+3654 q^4-6876 q^3+45792 q^2-10080 q+45360}{3780 q^8}, \\ 
b_3 &=\frac{179 q^8-1068 q^7+403 q^6+2184 q^5+20286 q^4-46656 q^3+195768 q^2-60480 q+181440}{40320 q^8}, \\ 
b_4 &=-\frac{-3 q^8+53 q^7-276 q^6+7 q^5+4158 q^4-9036 q^3+26676 q^2-10080 q+22680}{22680 q^8}, \\ 
b_5 &=\frac{3 q^8-100 q^7+1635 q^6-13160 q^5+58590 q^4-106560 q^3+236088 q^2-100800 q+181440}{1814400 q^8}.
\end{aligned}
\end{equation*}

When $N=6$:
\begin{equation*}
\tiny
\begin{aligned}
\\
a_0 &= \frac{51693 q^{10}-856 q^8+18721 q^6+1433070 q^4+10994040 q^2+9979200}{129600 q^{10}}, \\ 
a_1 &= -\frac{-73191 q^{10}+48926 q^9-22097 q^8+10306 q^7+35192 q^6-1980 q^5+2951028 q^4+1504080 q^3+22169520 q^2+1814400 q+19958400}{302400 q^{10}}, \\ 
a_2 &= \frac{13053 q^{10}+8834 q^9-21784 q^8+23242 q^7+5185 q^6+1476 q^5+1617966 q^4+1564560 q^3+11356920 q^2+1814400 q+9979200}{241920 q^{10}}, \\ 
a_3 &= -\frac{-4839 q^{10}-28638 q^9-6833 q^8+78966 q^7-69640 q^6+64908 q^5+3811572 q^4+4996080 q^3+23621040 q^2+5443200 q+19958400}{1088640 q^{10}}, \\
a_4 &= \frac{237 q^{10}+4372 q^9+24104 q^8+13892 q^7-93599 q^6+160200 q^5+2414574 q^4+3612960 q^3+12445560 q^2+3628800 q+9979200}{1814400 q^{10}}, \\ 
a_5 &= -\frac{-39 q^{10}-770 q^9-13937 q^8-111430 q^7-166408 q^6+1035540 q^5+6500340 q^4+9939600 q^3+26524080 q^2+9072000 q+19958400}{19958400 q^{10}}, \\ 
a_6 &= \frac{-3 q^{10}+198 q^9+2024 q^8+19998 q^7+239041 q^6+1324620 q^5+4548654 q^4+6629040 q^3+14259960 q^2+5443200 q+9979200}{119750400 q^{10}}, \\ 
b_1 &= \frac{73191 q^{10}+48926 q^9+22097 q^8+10306 q^7-35192 q^6-1980 q^5-2951028 q^4+1504080 q^3-22169520 q^2+1814400 q-19958400}{302400 q^{10}}, \\ 
b_2 &= \frac{13053 q^{10}-8834 q^9-21784 q^8-23242 q^7+5185 q^6-1476 q^5+1617966 q^4-1564560 q^3+11356920 q^2-1814400 q+9979200}{241920 q^{10}}, \\ 
b_3 &= \frac{4839 q^{10}-28638 q^9+6833 q^8+78966 q^7+69640 q^6+64908 q^5-3811572 q^4+4996080 q^3-23621040 q^2+5443200 q-19958400}{1088640 q^{10}}, \\ 
b_4 &= \frac{237 q^{10}-4372 q^9+24104 q^8-13892 q^7-93599 q^6-160200 q^5+2414574 q^4-3612960 q^3+12445560 q^2-3628800 q+9979200}{1814400 q^{10}}, \\ 
b_5 &= \frac{39 q^{10}-770 q^9+13937 q^8-111430 q^7+166408 q^6+1035540 q^5-6500340 q^4+9939600 q^3-26524080 q^2+9072000 q-19958400}{19958400 q^{10}}, \\ 
b_6 &= \frac{-3 q^{10}-198 q^9+2024 q^8-19998 q^7+239041 q^6-1324620 q^5+4548654 q^4-6629040 q^3+14259960 q^2-5443200 q+9979200}{119750400 q^{10}}.
\end{aligned}
\end{equation*}

\end{section}

\end{document}